\documentclass{amsart}
\usepackage{graphicx} 

\usepackage{hyperref}   
\usepackage[square,numbers]{natbib}     
\usepackage{amsthm}     
\usepackage{amsmath}    
\usepackage{amssymb}    
\usepackage{amsxtra}    
\usepackage{eucal}      
\usepackage[bbgreekl]{mathbbol}   
\usepackage{todonotes}
\usepackage{comment}
\usepackage{tikz}

\usetikzlibrary{matrix}
\usetikzlibrary{patterns, shapes}
\usetikzlibrary{arrows}
\usetikzlibrary{calc,3d}
\usetikzlibrary{decorations,decorations.pathmorphing, decorations.pathreplacing}
\usetikzlibrary{through}
\tikzset{every picture/.style={baseline=-.65ex}}
\tikzset{ext/.style={circle, draw,inner sep=1pt},int/.style={circle,draw,fill,inner sep=1pt},nil/.style={inner sep=1pt}}
\tikzset{exte/.style={circle, draw,inner sep=3pt},inte/.style={circle,draw,fill,inner sep=3pt}}
\tikzset{diagram/.style={matrix of math nodes, row sep=3em, column sep=2.5em, text height=1.5ex, text depth=0.25ex}}
\tikzset{diagram2/.style={matrix of math nodes, row sep=0.5em, column sep=0.5em, text height=1.5ex, text depth=0.25ex}}
\tikzset{every loop/.style={draw}}

\newcommand{\tadpole}{
\begin{tikzpicture}[baseline=-.65ex]
\node[int] (v) at (0,0) {};
\draw (v) edge[loop] (v);
\end{tikzpicture}
}

\usepackage{tikz-cd}

\theoremstyle{plain}

\newtheorem{thm}[subsection]{Theorem}

\newtheorem{prop}[subsection]{Proposition}
\newtheorem{lemma}[subsection]{Lemma}
\newtheorem{cor}[subsection]{Corollary}

\newcommand{\G}{\mathsf{G}}
\newcommand{\GC}{\mathsf{GC}}
\newcommand{\mF}{\mathcal{F}}
\DeclareMathOperator{\gr}{gr}
\DeclareMathOperator{\vdim}{dim}
\newcommand{\lp}{\text{-loop}}
\newcommand{\vt}{\text{-vert}}
\DeclareMathOperator{\ad}{ad}

\title{The triconnected Kontsevich graph complex}
\author{Thomas Willwacher}
\address{Department of Mathematics \\ ETH Zurich \\
R\"amistrasse 101 \\
8092 Zurich, Switzerland}
\email{thomas.willwacher@math.ethz.ch}

\begin{document}
\begin{abstract}
We show that a smaller version of the Kontsevich graph complex spanned by triconnected graphs is quasi-isomorphic to the full Kontsevich graph complex.
\end{abstract}

\maketitle

\section{Introduction}

The Kontsevich graph complex $\G_n$ is the space of linear combinations of pairs $(\gamma, o)$ consisting of a connected simple graph $\gamma$ with $\geq 3$-valent vertices and an orientation datum $o$, modulo natural relations.
It is equipped with a differential $d:\G_n\to \G_n$ by contracting edges,
\[
d(\gamma,o) = \sum_e (\gamma/e,o_e),
\]
see section \ref{sec:graphs} below for details.
The cohomology of $\G_n$, the Kontsevich graph cohomology, appears in a multitude of problems in algebra, topology and geometry.

A graph is called $k$-connected if it remains connected after deleting any subset of at most $k-1$ of its vertices, together with the adjacent edges.
A 1-connected graph is a connected graph, a 2-connected graph is also called biconnected, and a 3-connected graph is also called triconnected.
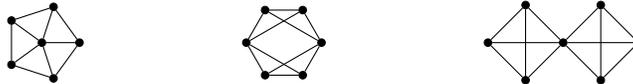
\begin{align*}
    &\begin{tikzpicture}
        \node[int] (v1) at (0:.5) {};
        \node[int] (v2) at (72:.5) {};
        \node[int] (v3) at (144:.5) {};
        \node[int] (v4) at (-144:.5) {};
        \node[int] (v5) at (-72:.5) {};
        \node[int] (c) at (0,0) {};
        \draw (c) edge (v1) edge (v2) edge (v3) edge (v4) edge (v5)
        (v1) edge (v2) edge (v5) (v3) edge (v2) edge (v4) (v4) edge (v5);
    \end{tikzpicture}
    &
    &\begin{tikzpicture}
        \node[int] (v1) at (0:.5) {};
        \node[int] (v2) at (60:.5) {};
        \node[int] (v3) at (120:.5) {};
        \node[int] (v4) at (180:.5) {};
        \node[int] (v5) at (-120:.5) {};
        \node[int] (v6) at (-60:.5) {};
        \draw (v1) edge (v3) (v2) edge (v4) (v4) edge (v6) (v5) edge (v1)
        (v1) edge (v2) edge (v6) (v3) edge (v2) edge (v4) (v5) edge (v6) edge (v4);
    \end{tikzpicture}
    &  &\begin{tikzpicture}
        \node[int] (v1) at (0,0) {};
        \node[int] (v2) at (1,0) {};
        \node[int] (v3) at (-1,0) {};
        \node[int] (v4) at (.5,.5) {};
        \node[int] (v5) at (.5,-.5) {};
        \node[int] (v6) at (-.5,.5) {};
        \node[int] (v7) at (-.5,-.5) {};
        \draw (v1) edge (v2) edge (v3) edge (v4) edge (v5) edge (v6) edge (v7)
        (v6) edge (v7) (v4) edge (v5)
        (v2) edge (v4) edge (v5)
        (v3) edge (v6) edge (v7);
    \end{tikzpicture}
\end{align*}
For example, the 5-wheel graph is triconnected, the graph in the middle is biconnected but not triconnected, and the graph on the right-hand side is neither.
The operation of contracting an edge can not increase the connectivity of a graph. Hence we have a natural sequence of quotients of the Kontsevich graph complex
\begin{equation}\label{equ:G chain}
\G_n \to \G_n^{bi} \to \G_n^{tri},
\end{equation}
where the biconnected Kontsevich graph complex 
\[
\G_n^{bi} = \G_n / I_2
\]
is the quotient by the subspace $I_2$ spanned by non-biconnected graphs and 
\[
\G_n^{tri} = \G_n / I_3
\]
is the quotient by the subspace $I_3$ spanned by the non-triconnected graphs. In fact, the natural Lie cobracket defined on $\G_n$ can at most reduce connectivity as well, so that $\G_n^{bi}$ and $\G_n^{tri}$ inherit the Lie coalgebra structure from $\G_n$, and the maps \eqref{equ:G chain} are morphisms of dg Lie coalgebras.

It is well-known (see \cite{CGV} or \cite[Appendix F]{grt}) that the projection 
\[
\G_n \to \G_n^{bi}
\]
is a quasi-isomorphism. In fact, in parts of the literature the biconnectivity condition is even imposed from the start and our $\G_n^{bi}$ is considered "the" Kontsevich graph comples.

Our main result is that the triconnected graph complex $\G_n^{tri}$ is also quasi-isomorphic to the Kontsevich graph complex $\G_n$.
\begin{thm}\label{thm:main}
    For every $n$ the projection 
    \[
    \G_n \to \G_n^{tri}
    \]
    is a quasi-isomorphism.
\end{thm}

Our result can be used to deduce improved lower bounds on the dimension of $H(\G_n)$ for even $n$ as we shall explain in section \ref{sec:applications} below.

    One may ask whether one can continue the chain \eqref{equ:G chain} to complexes of even more highly connected graphs. While it is straightforward to define these more highly connected complexes and extend  the chain, they will be far from quasi-isomorphic to the original Kontsevich graph complex. To see this, note that all vertices of a $k$-connected graph must have valence $\geq k$.
    Hence already the $4$-connected version of $\G_2$ is concentrated in positive cohomological degrees, while $H(\GC_2)$ is concentrated in non-positive degrees.
    Hence Theorem \ref{thm:main} is in this sense optimal -- its analogue for the more than triconnected graph complexes fails.

\subsection*{Acknowledgements}
The author is grateful for discussions with Peter Patzt, from which this project arose. 
This work has been partially supported by the NCCR Swissmap funded by the Swiss National Science Foundation, and the Horizon Europe Framework Programme CaLIGOLA (101086123).

\section{Background}
\subsection{The Kontsevich graph complex and its variants}
\label{sec:graphs}
Fix an integer $n$. Let $\gamma$ be a connected simple graph, that is, a combinatorial graph without multiple edges or self-edges. Then an $n$-orientation for $\gamma$ is the following data, depending on the parity of $n$:
\begin{itemize}
    \item For $n$ even, an $n$-orientation is an ordering of the set of edges of $\gamma$.
    \item For $n$ odd, an $n$-orientation is an ordering of the sets of half-edges and vertices $\gamma$.
\end{itemize}
We assign to a graph $\gamma$ with $r$ vertices and $k$ edges the cohomological degree $(n-1)k -n(r-1)$.
Then the Kontsevich graph complex $\G_n$ is the graded vector space spanned by pairs $(\gamma,o)$ consisting of a simple connected graph $\gamma$ whose vertices are $\geq 3$-valent and an $n$-orientation $o$ for $\gamma$, modulo the following two relations:
\begin{itemize}
    \item (Isomorphism) Let $\phi:\gamma\to \nu$ be a graph isomorphism. Let $o$ be an $n$-orientation for $\gamma$ and $\phi_*o$ the natural $n$-orientation of $\nu$ induced by $\phi$. Then we set 
    \[
    (\gamma,o) = (\nu,\phi_*(o)).
    \]
    \item (Orientation change) Suppose that $o'$ is another $n$-orientation of $\gamma$ obtained by applying some permutation $\pi$ to the ordering. Then we set 
    \[
    (\gamma,o')= \mathrm{sgn}(\sigma) (\gamma,o).
    \]
\end{itemize}
For even $n$ we write $o=e_1\wedge\cdots \wedge e_k$ for the $n$-orientation given by the ordering $e_1<\dots<e_k$ of the edge set $\{e_1,\dots,e_k\}$ of $\gamma$. Similarly, for $n$ odd we use the notation $o=v_1\wedge\cdots\wedge v_r\wedge h_1\wedge \cdots \wedge h_{2k}$ to define an $n$-orientation, where $\{v_1,\dots,v_r\}$ is the set of vertices and $\{h_1,\dots,h_{2k}\}$ the set of half-edges of $\gamma$.

The differential on $\G_n$ is given by edge contraction. For $n$ even we have
\[
d (\gamma,e_1\wedge\cdots\wedge e_k)=\sum_{j=1}^k (-1)^{j-1} (\gamma/e_j,e_1,\wedge\cdots\wedge \hat e_j\wedge \cdots \wedge e_k),
\]
where $\gamma/e_j$ is the graph obtained by contracting the edge $e_j$.
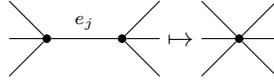
\[
\begin{tikzpicture}
\node[int] (v) at (0,0) {};
\node[int] (w) at (1,0) {};
\draw (v) edge +(-.5,.5) edge +(-.5,0) edge +(-.5,-.5) 
(w) edge +(.5,.5) edge +(.5,0) edge +(.5,-.5) 
(v) edge node[above] {$\scriptstyle e_j$} (w);
\end{tikzpicture}
\mapsto
\begin{tikzpicture}
\node[int] (v) at (0,0) {};
\draw (v) edge +(-.5,.5) edge +(-.5,0) edge +(-.5,-.5) 
(v) edge +(.5,.5) edge +(.5,0) edge +(.5,-.5) ;
\end{tikzpicture}
\]
If the contraction of the edge results in a graph that is not simple, we will implicitly set the corresponding term in the sum to zero.
We will write the above formula for short as 
\[
d (\gamma,o)= \sum_e (\gamma/e,o_e),
\]
with the sum running over the edges of $\gamma$.
The same formula applies also to the case of odd $n$ for a slightly altered definition of the orientation $o_e$ of $\gamma/e$, see \cite[eqn. (2)]{BW}.

We remark that there are several possible variations of the definition of the Kontsevich graph complex. For example, one may drop the simplicity condition on the graphs, and/or relax the trivalence condition on vertices. One can show that these variations do not affect the cohomology in loop orders $\geq 3$, see \cite{WZ, grt}.

\subsection{Graphs, connectivity, and SPQR trees}
Let $\gamma$ be a connected graph with at least 3 vertices. Then $\gamma$ is called \emph{biconnected} if it remains connected after removing any vertex, and \emph{triconnected} if it also remains connected after removing any pair of vertices, together with the adjacent edges.
We say that a pair of vertices $u\neq v$ of $\gamma$ is a \emph{separation pair} if $\gamma$ becomes disconnected after removing $u,v$, and their adjacent edges.

To any biconnected graph we can associate its SPQR tree. This is the following data:
\begin{itemize}
\item A tree $T$ whose nodes are of either of three types, S, P and R. 
\item No S or P nodes can have neighbors of the same type.
\item To each node $V$ of $T$ there is associated a \emph{skeleton graph} $T(V)$.
This is a graph on a subset $A_V$ of the set of vertices of $\gamma$. It can have two kinds of edges, virtual and real. The real edges must be a subset $E_V$ of the edges of $\gamma$.
The endpoints of the virtual edges form a separating pair of $\gamma$.
We will draw virtual edges as dotted and real edges as solid edges in drawings.
\item The skeleton $T(V)$ associated to an S node $V$ is a cycle graph with $\geq 3$ verices, for example 
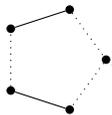
\[
\begin{tikzpicture}
\node[int] (v1) at (0:.7) {};
\node[int] (v2) at (72:.7) {};
\node[int] (v3) at (144:.7) {};
\node[int] (v4) at (-144:.7) {};
\node[int] (v5) at (-72:.7) {};
\draw (v2) edge[dotted] (v1) edge (v3) (v4) edge[dotted] (v3) edge (v5) (v1) edge[dotted] (v5);
\end{tikzpicture}
.
\]
\item The skeleton $T(V)$ associated to a P node $V$ is a graph with exactly 2 vertices connected by $\geq 3$ edges, for example 
\[
\begin{tikzpicture}
    \node[int] (v1) at (0,-.5) {};
\node[int] (v2) at (0,.5) {};
\draw (v1) edge[dotted, bend left] (v2) edge (v2) edge[dotted, bend right] (v2);
\end{tikzpicture}
.
\]
\item The skeleton $T(V)$ associated to an R node $V$ is a triconnected graph.
\item The half-edges of the tree $T$ incident to node $V$ are in 1-1 correspondence to the virtual edges of $T(V)$. This association must be such that the respective endpoints must correspond to the same vertices of $\gamma$.
\item The union of the vertices of the skeleta is $\gamma$.
The edges of $\gamma$ are the disjoint union of the sets of real edges $E_V$.
\end{itemize}

\begin{figure}
\begin{tikzpicture}[scale=.5]
    \node[int] (s1) at (0,1) {};
    \node[int] (s2) at (1,0) {};
    \node[int] (s3) at (0,-1) {};
    \node[int] (s4) at (-1,0) {};
    \node[int] (a1) at (-2,1) {};
    \node[int] (a2) at (-1,2) {};
    \node[int] (b1) at (2,1) {};
    \node[int] (b2) at (1,2) {};
    \node[int] (c1) at (-2,-1) {};
    \node[int] (c2) at (-1,-2) {};
    \draw (a1) edge (a2) edge (s1) edge (s4)
          (a2) edge (s1) edge (s4)
          (b1) edge (b2) edge (s1) edge (s2)
          (b2) edge (s1) edge (s2)
          (c1) edge (c2) edge (s3) edge (s4)
          (c2) edge (s3) edge (s4)
          (s3) edge (s4) edge (s2);
\end{tikzpicture}
$\quad \leftrightarrow\quad$
\begin{tikzpicture}[scale=.5]
\begin{scope}[scale=.7, xshift=-1.5cm, yshift=1.5cm]
    \node[int] (a1) at (-2,1) {};
    \node[int] (a2) at (-1,2) {};
    \node[int] (a3) at (0,1) {};
    \node[int] (a4) at (-1,0) {};
    \draw (a1) edge (a2) edge (a3) edge (a4)
        (a2) edge (a3) edge (a4)
        (a3) edge [dashed] (a4);
    \node[draw, circle, inner sep=.35cm,label={R}] (R1) at (-1,1) {};
\end{scope}
\begin{scope}[scale=.7, xshift= 3.5cm, yshift=1.5cm]
    \node[int] (b1) at (-2,1) {};
    \node[int] (b2) at (-1,2) {};
    \node[int] (b3) at (0,1) {};
    \node[int] (b4) at (-1,0) {};
    \draw (b1) edge (b2) edge (b3) edge[dashed] (b4)
        (b2) edge (b3) edge (b4)
        (b3) edge  (b4);
    \node[draw, circle, inner sep=.35cm,label={R}] (R2) at (-1,1) {};
\end{scope}
\begin{scope}[scale=.7,xshift= -3.5cm, yshift=-5.5cm]
    \node[int] (c1) at (-2,1) {};
    \node[int] (c2) at (-1,2) {};
    \node[int] (c3) at (0,1) {};
    \node[int] (c4) at (-1,0) {};
    \draw (c1) edge (c2) edge (c3) edge (c4)
        (c2) edge[dashed] (c3) edge (c4)
        (c3) edge  (c4);
    \node[draw, circle, inner sep=.35cm,label=-90:{R}] (R3) at (-1,1) {};
\end{scope}
\begin{scope}[scale=.7]
    \node[int] (s1) at (0,1) {};
    \node[int] (s2) at (1,0) {};
    \node[int] (s3) at (0,-1) {};
    \node[int] (s4) at (-1,0) {};
    \draw (s1) edge[dashed] (s2) edge[dashed] (s4) 
            (s3) edge (s2) edge[dashed] (s4);
    \node[draw, circle, inner sep=.35cm,label=-45:{S}] (S) at (0,0) {};
\end{scope}
\begin{scope}[yshift=-1.5cm, xshift=-1.5cm, scale=.7]
\node[int] (p1) at (-.5,.5) {};
\node[int] (p2) at (.5,-.5) {};
\draw (p1) edge[bend left, dashed] (p2) 
           edge[bend right, dashed] (p2) 
           edge (p2);
\node[draw, circle, inner sep=.25cm,label=-45:{P}] (P) at (0,0) {};
\end{scope}
\draw (S) edge (R1) edge (R2) edge (P) (P) edge (R3);
\end{tikzpicture}
    \caption{\label{fig:SPQR example} A biconnected graph (left) and its associated SPQR tree (right). The graph inscribed into a vertex of the SPQR tree is the skeleton of that vertex. The real edges are drawn as solid edges, and the virtual edges as dashed edges.}
\end{figure}
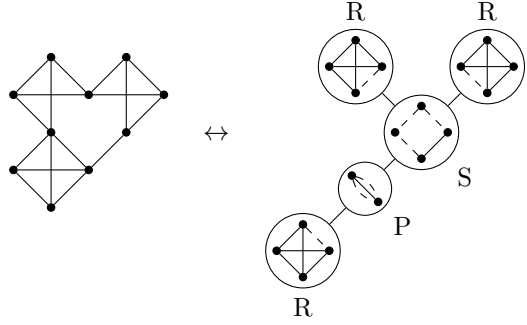

The graph $\gamma$ may be recovered from the SPQR tree by gluing together the matching virtual edges, and then removing those edges.
The SPQR tree of a biconnected graph is unique. It may be found by iteratively cutting the graph at separation pairs, and finally fusing adjacent S nodes and adjacent P nodes in the resulting tree.
We refer the reader to \cite{DBT} for more details on SPQR trees. For an example, see Figure \ref{fig:SPQR example}.

\begin{lemma}\label{lem:spqr leaves}
    Let $\gamma$ be a simple graph with at least trivalent vertices. Then every leaf of the SPQR tree $T$ of $\gamma$ is an R node. In particular, $T$ always has at least one R node.
\end{lemma}
\begin{proof}
    The leaves of $T$ are the nodes whose skeleta have only a single virtual edge.
    Now if (the skeleton of) a P node has only a single virtual edge then it needs to have at least two real edges. Hence the graph would not be not simple.
    If an S node $V$ only has a single virtual edge, then it has to have two consecutive real edges in the cycle $T(V)$, and hence the graph $\gamma$ has a bivalent vertex.  
\end{proof}

\section{Proof of Theorem \ref{thm:main}}
Since we already know that $\G_n \to \G_n^{bi}$ is a quasi-isomorphism it is enough to show that the projection 
\begin{equation}\label{equ:tbs}
\G_n^{bi}\to \G_n^{tri}
\end{equation}
is a quasi-isomorphism as well.
Consider the ascending exhaustive filtration
\begin{align*}
    0=\mF_{-1}\G_n^{bi} \subset \mF_0\G_n^{bi}
    \subset \mF_1\G_n^{bi}\subset \cdots
\end{align*} 
such that $\mF_p\G_n^{bi}$ is spanned by graphs $\gamma$ whose SPQR tree has $\leq p$ real edges in (the skeleta of) R-nodes. 
We will see below that 
\begin{equation}\label{equ:tbs compat}
d \mF_p \G_n^{bi} \subset \mF_p \G_n^{bi},
\end{equation}
so that each $\mF_p \G_n^{bi}$ defines a subcomplex of $\G_n^{bi}$.

The same definition also yields subspaces $\mF_p\G_n^{tri}\subset \G_n^{tri}$ that together form a filtration on $\G_n^{tri}$. Note that this latter filtration is simply the filtration by the number of edges since the SPQR tree of a triconnected graph consists of a single R node, and the skeleton of that node is the whole graph.
In particular, is it hence clear that 
\begin{equation}\label{equ:tri compat}
d \mF_p \G_n^{tri} \subset \mF_{p-1} \G_n^{tri}
\end{equation}

By \eqref{equ:tri compat} the differential on the associated graded complex $\gr \G_n^{tri}$ is zero.
Let us analyze the differential on $\gr \G_n^{bi}$.
To this end consider a biconnected graph $\gamma$ with $\geq 3$-valent vertices. Let $T$ be the SPQR tree of $\gamma$ and fix an edge $e$ of $\gamma$.
Let $V$ be the unique node of $T$ whose skeleton $T(V)$ contains $e$, i.e., $e$ is a real edge of $T(V)$. We then study the SPQR tree of the graph $\gamma/e$.\footnote{We note that SPQR tree updates upon edge contraction have been studied in the combinatorics literature \cite{KRGIH}.}
\begin{itemize}
    \item If $V$ is a P-node, then the endpoints of $e$ are the two vertices $u,v$ of the skeleton $T(V)$. The pair $(u,v)$ is a separation pair and hence contracting $e$ yields a non-biconnected graph. (The vertices $u,v$ are merged into a cut vertex of $\gamma/e$.)
    \item Assume $V$ is an S node, whose skeleton $T(V)$ is a $k$-cycle.
    \begin{itemize}
        \item If $k\geq 4$, then the SPQR tree of $\gamma/e$ is the same as $T$, but with the skeleton associated to $V$ set to $T(V)/e$.
        \item If $k=3$, then $T(V)$ must have the form 
        
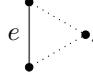
\[
            \begin{tikzpicture}
     \node[int] (v1) at (0:.5) {};   
     \node[int] (v2) at (120:.5) {};   
     \node[int] (v3) at (-120:.5) {};   
    \draw (v1) edge [dotted] (v2) edge [dotted] (v3) (v2) edge node[left] {$e$} (v3);
    \end{tikzpicture}.
        \]
        Note in particular that due to the trivalence condition on $\gamma$ the two edges adjacent to $e$ in the loop $T(V)$ must be virtual.
        Hence $V$ has two neighbors $V'$, $V''$ in $T$. If both $V'$ and $V''$ are $P$-nodes then the SPQR tree of $\gamma/e$ is obtained by removing $V$ from $T$ and merging $V'$ and $V''$ to one P node.
        Otherwise, the SPQR tree of $\gamma/e$ is obtained by removing $V$ and just connecting $V$ to $V''$.
    \end{itemize}
    
    \item If $V$ is an R node, then let $T'$ be the SPQR tree of $T(V)/e$. Then the SPQR tree of $\gamma/e$ is obtained from $T$ by inserting $T'$ in the place of node $V$, cf. \cite[Algorithm 2]{KRGIH}. The details of this gluing prodedure are not important here. We merely note that the resulting SPQR tree of $\gamma/e$ has at least one real edge less in R nodes than that of $\gamma$.
\end{itemize}

Note that \eqref{equ:tbs compat} follows directly from the above discussion. 
Furthermore, one sees that the differential of the associated graded $\gr \G_n^{bi}$ is the map 
\[
d' (\gamma,o)= \sum_e' (\gamma/e,o_e),
\]
where the sum is restricted to run over edges that are in S nodes in the SPQR tree of $\gamma$.

Now by elementary homological algebra results we know that \eqref{equ:tbs} is a quasi-isomorphism if we can show the associated graded morphism 
\begin{equation}\label{equ:tbs2}
\gr \G_n^{bi}\to \gr \G_n^{tri}
\end{equation}
is a quasi-isomorphism. 
This in turn is equivalent to showing that 
\[
H(\gr K_n) =0
\]
with 
\[
K_n = \ker (\G_n^{bi}\to \G_n^{tri}) \subset \G_n^{bi}.
\]
Concretely, $K_n$ is spanned by graphs whose SPQR tree has at least two nodes.
By Lemma \ref{lem:spqr leaves} we know that at least one\footnote{In fact at least two, but one is sufficient for our argument.} of these is an R leave.
The differential on $\gr K_n$ does not change the number of R nodes, nor their skeleta. We may hence consider separately the subcomplex $L_n^{(k)}\subset \gr K_n$ spanned by graphs with $k$ many R leaves, so that 
\[
\gr K_n = \bigoplus_{k\geq 1} L_n^{(k)}.
\]
Again since the differential $d'$ does not affect the R nodes we may temporarily pass to a bigger complex $L_n^k$ spanned by graphs equipped with the extra datum of an ordering on the $k$ R leaves of their SPQR tree. Then 
\[
L_n^{(k)} = (L_n^k)_{S_k}
\]
and it suffices to check that $H(L_n^{k})=0$.
To show this in turn let us define a homotopy
\[
h: L_n^k \to L_n^k
\]
that acts on a graph $\gamma$ with SPQR tree $T$ as follows. Let $V$ be the first R leaf, $e\in T(V)$ be its unique virtual edge and $W$ the unique neighbor of $V$.
\begin{itemize}
    \item If $W$ is an $S$-node, then we define $h\gamma$ by summing over all ways of inserting a real edge between any pair of consecutive virtual edges of $T(W)$.
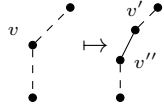
    \[
    \begin{tikzpicture}
        \node[int,label=170:{$\scriptstyle v$}] (a) at (0,0) {};
        \node[int] (b) at (0,-.7) {};
        \node[int] (c) at (.5,.5) {};
        \draw (a) edge[dashed] (b) edge [dashed] (c);
    \end{tikzpicture}
    \mapsto 
        \begin{tikzpicture}
        \node[int,label=0:{$\scriptstyle v''$}] (a) at (0,-.2) {};
        \node[int,label=90:{$\scriptstyle v'$}] (aa) at (.2,.2) {};
        \node[int] (b) at (0,-.7) {};
        \node[int] (c) at (.5,.5) {};
        \draw (a) edge[dashed] (b) (aa) edge [dashed] (c) (a) edge (aa);
    \end{tikzpicture}
    \]
    The orientation of the newly formed graph for $n$ even is such that the new edge becomes the first in the ordering of edges. For odd $n$ let $v$ be the old vertex, which we assume to be the first in the ordering, $v',v''$ the new vertices and $h',h''$ the new half-edges, with $h'$ adjacent to $v'$. Then the orientation of the new graph can be taken to be $h'\wedge h''\wedge v'\wedge v''\wedge (\dots)$, with all other vertices and half-edges retaining their position in the ordering.
    \item If $W$ is either a P or an R node, then we insert a new S node $U$ between $V$ and $W$ with skeleton
    
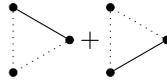
\[
    \begin{tikzpicture}
     \node[int] (v1) at (0:.5) {};   
     \node[int] (v2) at (120:.5) {};   
     \node[int] (v3) at (-120:.5) {};   
    \draw (v1) edge (v2) edge [dotted] (v3) (v2) edge[dotted] (v3);
    \end{tikzpicture}
    +
        \begin{tikzpicture}
     \node[int] (v1) at (0:.5) {};   
     \node[int] (v2) at (120:.5) {};   
     \node[int] (v3) at (-120:.5) {};   
    \draw (v1) edge[dotted] (v2) edge (v3) (v2) edge [dotted] (v3);
    \end{tikzpicture}
    \]
    In other words, the operation is the same as in the first case, but after we formally create an S node neighbor of $V$ with skeleton 
    \begin{equation}\label{equ:triv skeleton}
    \begin{tikzpicture}
        \node[int] (v) at (0,.5) {};
        \node[int] (w) at (0,-.5) {};
        \draw (v) edge[dotted, bend left] (w) edge[dotted, bend right] (w);
    \end{tikzpicture}.
    \end{equation}
\end{itemize}

We next claim that for any graph $\gamma$ in $K_n^k$ we have that
\begin{equation}\label{equ:homotopy}
d'h(\gamma) + h(d'\gamma) = 
N_\gamma \gamma,
\end{equation}
where $N_\gamma$ is the number of virtual edges in the S node $W$ adjacent to the first R leave $V$, or $N_\gamma=2$ if $W$ is not an S node.
In fact, to unify the discussion, we will always consider the case that $W$ is an S node. If it is not, we formally insert an S node with skeleton \eqref{equ:triv skeleton} and make that node the new $W$.
Then we split $d'=d_W+d''$, with $d_W$ the terms of $d'$ contracting edges in the S node $W$, and $d''$ the terms contracting edges in S nodes elsewhere in the graph. Then we clearly have 
\[
d''h(\gamma) + h(d''\gamma) = 0
\]
and 
\[
d_Wh(\gamma) + h(d_W\gamma) = N_\gamma \gamma
\]
so that \eqref{equ:homotopy} follows.
Finally, we note for all nonzero $\gamma$ we have that $N_\gamma>0$, and hence by \cite[Lemma 7]{FFW} we conclude that $H(L_n^k)=0$, and hence $H(\gr K_n)=0$ as desired. 

\hfill\qed

\section{Applications}\label{sec:applications}
\subsection{Lower dimensionality}
For computational purposes it is obviously desirable to have as small as possible complexes computing the Kontsevich graph cohomology.
In this regard the complex $\G_n^{tri}$ introduced in this paper is (slightly) preferable over $\G_n^{bi}$ and $\G_n$. Table \ref{table1} lists the dimensions of the three complexes at loop order 10, which is currently the maximal loop order in which the cohomology can be numerically computed. We see that our complex is roughly 15-20\% smaller than $\G_n^{bi}$ and $\G_n$.

\begin{figure}
\begin{tabular}{|c|c|c|c|}
\hline
$k$ & $\vdim \G_0^{k\vt,10\lp}$ & $\vdim \G_0^{bi, k\vt,10\lp}$ & $\vdim \G_0^{tri, k\vt,10\lp}$ \\
\hline
 6 & 1 & 1 (100\%) & 1 (100\%) \\
 7 & 4 & 4 (100\%) & 4 (100\%) \\
 8 & 291 & 291 (100\%) & 284 (97\%) \\
 9 & 5849 & 5846 (99\%) & 5461 (93\%) \\
 10 & 50863 & 50766 (99\%) & 45662 (89\%) \\
 11 & 243743 & 242778 (99\%) & 212426 (87\%) \\
 12 & 716800 & 712762 (99\%) & 611249 (85\%) \\
 13 & 1364145 & 1354993 (99\%) & 1144706 (83\%) \\
 14 & 1716410 & 1703974 (99\%) & 1423137 (82\%) \\
 15 & 1421127 & 1410740 (99\%) & 1167862 (82\%) \\
 16 & 746494 & 741356 (99\%) & 609526 (81\%) \\
 17 & 226296 & 224974 (99\%) & 183983 (81\%) \\
 18 & 30307 & 30187 (99\%) & 24585 (81\%) \\
 \hline
\end{tabular}
    \caption{\label{table1}Table of the dimensions of the 10-loop part of the three complexes $\G_0$, $\G_0^{bi}$ and $\G_0^{bi}$ for various numbers of vertices $k$. The percentage is the fraction of the size of $\G_0$.}
\end{figure}

\subsection{Spectral sequence and lower bounds on cohomology for even $n$}
Let us consider the dual graph complexes (and dg Lie algebras)
\begin{align*}
 \GC_n&:= \G_n^*
 &
 \GC_n^{bi}&:= (\G_n^{bi})^* 
 &
 \GC_n^{tri}&:= (\G_n^{tri})^*.
\end{align*}
Elements of these complexes can be seen as series of graphs with an orientation as before. However, dual differential $\delta_0$ acts through vertex splitting instead of edge contraction. The subcomplexes (and dg Lie subalgebras)
\[
 \GC_n^{tri} \subset \GC_n^{bi} \subset \GC_n
\]
consist of series of triconnected (resp. biconnected) graphs.

Let us focus on the case of even $n$. Specifically, we choose $n=0$.
Following \cite{KWZ1}, let $\nabla:\GC_0\to \GC_0$ be the degree +1 operation of adding one edge to a graph. In terms of the Lie bracket, this can be written as $\nabla \Gamma=[\tadpole,\Gamma]$. 
The operator $\nabla$ increases the loop order of graphs by $+1$.
Then one has:

\begin{thm}[{\cite[Section 3 up to Corollary 4]{KWZ1}}]\label{thm:KWZ}
    We have that $(\delta_0+\nabla)^2=0$ and
    \[
    H(\GC_0,\delta_0+\nabla) = \bigoplus_{k\geq 1} \mathbb{Q} \sigma_{2k+1}
    \]
    with $\sigma_{2k+1}$ a cohomology class of loop order $\leq 2k+1$ and degree (number of edges) $4k+2$.
\end{thm}

Now there is an obvious convergent spectral sequence on $(\GC_0,\delta_0+\nabla)$ arising from the filtration by loop order such that 
\[
E_1 = H(\GC_0, \delta_0) \Rightarrow H(\GC_0,\delta_0+\nabla).
\]
We call the above spectral sequence the \emph{loop order spectral sequence}.
The idea of \cite{KWZ1} is to exploit this spectral sequence to obtain information on the (unknown) Kontsevich graph cohomology $H(\GC_0, \delta_0)$ from the known cohomology $H(\GC_0,\delta_0+\nabla)$.
As a corollary of our main result Theorem \ref{thm:main} we may show:

\begin{cor}\label{cor:sseq}
    Let $d_k$ be the differential on the $E_k$-page of the loop order spectral sequence. (For example, $d_0=\delta_0$ and $d_1$ is induced by $\nabla$.) Then we have 
    \[
    d_1=d_3=d_5=0.
    \]
\end{cor}

\subsubsection{Dimension bounds}
Before we continue with the proof of Corollary \ref{cor:sseq}, let us briefly discuss the dimension bounds on the graph cohomology $H(\GC_0,\delta_0)$ that it implies.
To better compare our results with the literature, let us also discuss in parallel the cohomology $H(\GC_2,\delta_0)$ that is isomorphic up to degree shifts,
\[
H^k(\GC_2^{g\lp},\delta_0)
=
H^{k+2g}(\GC_0^{g\lp},\delta_0).
\]
\newcommand{\FL}{\mathfrak{g}}
\newcommand{\grt}{\mathfrak{grt}}
\newcommand{\fh}{\mathfrak{h}}
\newcommand{\fk}{\mathfrak{k}}
\newcommand{\mL}{\mathcal{L}}
One has the following results from the literature:
\begin{itemize}
    \item $H^k(\GC_2^{g\lp})=0$ for $k<0$ or $k>g-3$, or equivalently 
    \begin{equation}
        \label{equ:GC0 van}
        H^k(\GC_0^{g\lp})=0 \quad\quad \text{for $k<2g$ or $k>3g-3$.}
    \end{equation}
    \item $H^0(\GC_2)\cong \grt_1$ can be identified with the Grothendieck-Teichmüller Lie algebra \cite{grt}. Furthermore, by the work of F. Brown \cite{Brown} there is an injective map 
    \[
    \mathrm{FreeLie}(\sigma_3',\sigma_5',\dots) \to \grt_1, 
    \]
    from a free Lie algebra with generators $\sigma_{2k+1}'$ of cohomological degree 0 and loop order $2k+1$.
    For $\GC_0$ this implies that there is an injective map 
    \[
    \FL:=\mathrm{FreeLie}(\sigma_3,\sigma_5,\dots) \to H(\GC_0),
    \]
    where the generators $\sigma_{2k+1}$ are of cohomological degree (number of edges) $4k+2$ and loop order $2k+1$.
\end{itemize}
The generators $\sigma_{2k+1}$ above correspond to the $\sigma_{2k+1}$ of Theorem \ref{thm:KWZ}. 
In particular, picking some $(\delta_0+\nabla)$-cocycles representing the $\sigma_{2k+1}$ we obtain a morphism 
\[
\FL \to H(\GC_0, \delta_0+\nabla).
\]
By Theorem \ref{thm:KWZ} we have $H(\GC_0, \delta_0+\nabla)\cong \FL/[\FL,\FL]$, and the map above agrees with the natural projection $\FL\to \FL/[\FL,\FL]$ under this identification.
In particular the restriction to the commutator ideal $\fh=[\FL,\FL]$ 
\[
\fh = [\FL,\FL] \to H(\GC_0, \delta_0+\nabla)
\]
is the zero map. 
Let $E_k$ be the $k$-th page of the loop order spectral sequence from above. We then define for $k=1,2,\dots$
\[
\mL_k \fh := \ker(\fh \to E_k).
\]
We obviously have 
\[
0= \mL_1 \fh \subset \mL_2 \fh \subset \cdots 
\]
and 
\[
\cup_k \mL_k \fh = \fh.
\]

The elements of the associated graded 
\[
\fk_{k} := \mL_{k+1} \fh/\mL_{k} \fh
\]
can be understood as cohomology classes in the image of $\fh$ in $H(\GC_0,\delta_0)$ that survive in the loop order spectral sequence up to and including the $k$-th page, and are then killed by the differential $d_k$ and become zero on the $k+1$st page. Our Corollary \ref{cor:sseq} hence implies that 
\begin{equation}\label{equ:fk vani}
\fk_{1} = \fk_{3} = \fk_{5} = 0.
\end{equation}
We note several further facts:

\begin{itemize}
\item We have the lower dimension bound 
\begin{equation}\label{equ:lower bound}
\vdim H^{2g-1}(\GC_0^{g-k\lp},\delta_0) \geq \vdim \fk_{k}^{g\lp} 
\end{equation}
since the class on the $E_k$ page that kills our classes in $\fh$ have to be present already on the $E_1$ page, in the appropriate degree and loop order.
    \item By compatibility of our spectral sequence with the Lie brackets one has 
    \begin{equation}\label{equ:Lie compat}
    [\mL_k \fh, \FL] \subset \mL_{k} \fh.
    \end{equation}
    \item From \cite[Proposition 5]{KWZ1} one obtains that 
    \begin{equation}\label{equ:sigma fil}
    [\sigma_{2k+1},\FL] \subset \mL_{2k+1}\fh.
    \end{equation}
\end{itemize}

Also note that
\begin{equation}\label{equ:mL sum}
\vdim \mL_k \fh = \sum_{j=1}^{k-1} \vdim\fk_j.
\end{equation}

By numerical computation one can determine $H(\GC_0^{g\lp},\delta_0)$ for $g\leq 10$ \cite{BW}, see Figure \ref{fig:GC cohom}.

\begin{figure}
\begin{tabular}{c|ccccccccc}
    $k,g$ & 3 & 4 & 5 & 6 & 7 & 8 & 9 & 10 \\
    \hline
    0 & 1 & 0 & 1 & 0 & 1 & 1 & 1 & 1  \\
    \cline{2-2}
    1 & 0 & \multicolumn{1}{|c}{0} & 0 & 0 & 0 & 0 & 0 & 0  \\
    \cline{3-3}
    2 & 0 & 0 & \multicolumn{1}{|c}{0} & 0 & 0 & 0 & 0 & 0  \\
    \cline{4-4}
    3 & 0 & 0 & 0 & \multicolumn{1}{|c}{1} & 0 & 1 & 1 & 2  \\
    \cline{5-5}
    4 & 0 & 0 & 0 & 0 & \multicolumn{1}{|c}{0} & 0 & 0 & 0  \\
    \cline{6-6}
    5 & 0 & 0 & 0 & 0 & 0 & \multicolumn{1}{|c}{0} & 0 & 0  \\
    \cline{7-7}
    6 & 0 & 0 & 0 & 0 & 0 & 0 & \multicolumn{1}{|c}{0} & 0  \\
    \cline{8-8}
    7 & 0 & 0 & 0 & 0 & 0 & 0 & 0 & \multicolumn{1}{|c}{1}      
\end{tabular}
    \caption{\label{fig:GC cohom} The table shows $\vdim H^k(\GC_2^{g\lp})$ as computed in \cite{BW}. There is no cohomology in negative degrees for any $g$, and no cohomology in degrees $k>g-3$, as indicated by the diagonal line. }
\end{figure}

We furthermore have that 
\[
\mL_{3}\fh^{12\lp} = \fh^{12\lp}.
\]
Indeed, otherwise the corresponding classes would be canceled on the $E_4$ page of our spectral sequence or later. Hence there would need to be cohomology in $H(\GC_0)$ in loop orders $\leq 8$ and degree $k=2\cdot 12-1=23$. But this is impossible by \eqref{equ:GC0 van}. 
In particular, we have $[\sigma_5,\sigma_7]\in \mL_{2}\fh$.
From \eqref{equ:Lie compat} and \eqref{equ:sigma fil} we then have
\[
\langle [\sigma_5,\sigma_7],\sigma_3\rangle \cap \fh
\subset \mL_{3}\fh,
\]
where $\langle X \rangle$ denotes the ideal of $\FL$ generated by $X$.
Similarly, we obtain 
\begin{align*}
    \langle\sigma_3,\sigma_5\rangle \cap \fh
& \subset \mL_{5}\fh 
\\
    \langle\sigma_3,\sigma_5,\sigma_7\rangle \cap \fh
& \subset \mL_{7}\fh. 
\end{align*}
From \eqref{equ:fk vani}, \eqref{equ:lower bound} and \eqref{equ:mL sum} we then obtain the following lower dimension bounds on the Kontsevich graph cohomology.
\begin{cor}
\begin{itemize}
\item (Strengthening of \cite[Proposition 2.4]{PayneWillwacher})
    Let $a_g$ be the dimension of the $g$-loop part of $\langle [\sigma_5,\sigma_7],\sigma_3\rangle \cap \fh$. Then
    \[
    \vdim H^3(\GC_2^{g-2\lp},\delta_0) \geq a_{g}.
    \]
\item Let $b_g$ be the dimension of the $g$-loop part of $\langle\sigma_3,\sigma_5\rangle \cap \fh$. Then
    \[
    \vdim H^3(\GC_2^{g-2\lp},\delta_0) 
    +
    \vdim H^7(\GC_2^{g-4\lp},\delta_0) 
    \geq b_{g}.
    \]
\item Let $c_g$ be the dimension of the $g$-loop part of $\langle\sigma_3,\sigma_5,\sigma_7\rangle \cap \fh$. Then
    \[
    \vdim H^3(\GC_2^{g-2\lp},\delta_0) 
    +
    \vdim H^7(\GC_2^{g-4\lp},\delta_0) 
        +
    \vdim H^{11}(\GC_2^{g-6\lp},\delta_0)
    \geq c_{g}.
    \]
    \end{itemize}
\end{cor}

\subsubsection{Proof of Corollary \ref{cor:sseq}}

We shall temporarily need to work over a larger complex 
\[
\GC_0^{disc} \supset \GC_0
\]
that is defined in the same way as $\GC_0$ except that we only require graphs to be simple and $\geq 3$-valent, but we allow non-connected graphs. The operations $\delta_0$ and $\nabla$ remain well-defined over the larger complex $\GC_0^{disc}$. There is an obviuos projection 
\[
\pi:\GC_0^{disc} \to \GC_0 
\]
sending all disconnected graphs to zero, which satisfies $\delta_0\pi=\pi\delta_0$, since $\delta_0$ cannot alter the number of connected components. However, we have $\nabla\pi\neq\pi\nabla$.
Next, we need to recall the following result of \v Zivkovi\'c:
\begin{prop}[{\cite[section 3.1]{Z}}]\label{prop:Z}
    There is an operation $D:\GC_0^{disc} \to \GC_0^{disc}$ of cohomological degree 0 with the following properties:
    \begin{itemize}
        \item For $\gamma$ a $k$-connected graph of loop order $g$ we have that $D\gamma$ is a linear combination of $k-1$-connected graphs of loop order $g+1$. 
        \item We have the following identities:
        \begin{align*}
            \delta_0 D+D\delta_0 &= \nabla &
            D^2 &= 0 &
            \nabla D + D\nabla &= 0. 
        \end{align*}
    \end{itemize}
\end{prop}
Concretely, the operation $D$ is defined by removing one vertex and connecting the incident edges to other vertices in all possible ways, see \cite[]{Z}. (The precise combinatorial form of $D$ will not be relevant here.)

We also define the operation 
\[
\bar D := \pi D \colon \GC_0  \to \GC_0.
\]
It follows that 
\[
\delta_0 \bar D +\bar D\delta_0 = \pi(\delta_0D+D\delta_0)= \pi \nabla = \nabla.
\]
Furthermore, note that for biconnected $\gamma\in \GC_n^{bi}$ we have that $D\gamma$ is connected and hence
\begin{align}\label{equ:bar D prop}
\bar D^2 \gamma &= \pi D^2\gamma =0
&
(\nabla \bar D +\bar D \nabla) \gamma &= \nabla D\gamma + D \nabla \gamma =0.
\end{align}

Using $\bar D$ we can produce an isomorphism of (loop-order-)filtered complexes 
\begin{equation}\label{equ:cx iso}
\exp(\bar D) \colon (\GC_0,\delta_0+\nabla) \to  (\GC_0,\tilde \delta)
\end{equation}
with 
\[
\tilde \delta = e^{\bar D} (\delta_0+\nabla)e^{-\bar D}
=
\sum_{j\geq 0} \frac{1}{j!} \ad_{\bar D}^j (\delta_0+\nabla)
=
\delta_0+\delta_2+\delta_3+\cdots
\]
and 
\[
\delta_k= \frac{1}{k!} \ad_{\bar D}^k \delta_0 + \frac{1}{(k-1)!} \ad_{\bar D}^{k-1} \nabla = \frac{k-1}{k!} \ad_{\bar D}^{k-1} \nabla
\]
the part of the differential that increases the loop order by $k$.
From \eqref{equ:bar D prop} we then immediately obtain that 
\begin{equation}\label{equ:delta van}
\delta_k \gamma =0 \quad\quad  \text{ for all $\gamma\in \GC_n^{bi}$ and $k\geq 3$.}
\end{equation}
We will furthermore use the following elementary Lemma.
\begin{lemma}\label{lem:aux}
Let $(V,d)$ be a cochain complex and $U\subset V$ a subcomplex, such that the inclusion is a quasi-isomorphism. Then for any $x\in V$, $y\in U$ solving the equation 
\[
dx =y
\]
there is a $z\in V$ such that $x+dz\in U$.
\end{lemma}
\begin{proof}
    It is clear that $dy=0$, and that the cohomology class represented by $y$ in $H(V)$ is zero. Hence, since the inclusion $U\subset V$ is a quasi-isomorphism, there must be some $u\in U$ such that $du=y$. Hence $d(x-u)=0$. The cohomology class represented by $x-u$ in $H(V)$ must be the image of some cohomology class in $H(U)$, again since $U\subset V$ is a quasi-isomorphism. Hence there is some $z\in V$ such that $x-u+dz\in U$. Hence $x+dz\in U$ as well, as desired.
\end{proof}

\begin{proof}[Proof of Corollary \ref{cor:sseq}]
Since \eqref{equ:cx iso} is an isomorphism of filtered complexes both sides have isomorphic spectral sequences. We may hence consider the spectral sequence associated to the loop order filtration on $(\GC_0,\tilde \delta)$.
First, since the differential $\tilde \delta=\delta_0+\delta_2+\cdots$ has no part that increases the loop order by one, it is clear that the differential $d_1$ on the $E_1$-page of the spectral sequence vanishes, $d_1=0$.
We next show the differential on the $E_5$-page also vanishes, $d_5=0$.
This statemenent means that if we are given an element $X=X_0+X_1+X_2+X_3+X_4\in \GC_0$, with $X_j\in \GC_0^{g+j\lp}$ such that $\tilde \delta X \in \GC_0^{\geq g+5\lp}$, then there is an extension $X'=X_0+X_1'+X_2'+X_3'+X_4'+X_5'$ such that $\tilde \delta X' \in \GC_0^{\geq g+6\lp}$.
The closedness condition for $X$ translates into the following loop order components:
\begin{align}
    \delta_0 X_0 &= 0 
    &
    \delta_0 X_1 &= 0
    \\
    \delta_2 X_0 + \delta_0 X_2 &=0
    &
    \delta_3 X_0 +  \delta_2 X_1 + \delta_0 X_3 &=0
    \\
    \delta_4 X_0 +  \delta_3 X_1 + \delta_2 X_2 + \delta_0 X_4 &=0.\label{equ:g4}
\end{align}
We may also assume from the start that $X_0\in \GC_0^{tri}$ is triconnected, since by Theorem \ref{thm:main} we may otherwise change $X$ by an exact element to ensure this condition. In fact, by the same reasoning we may also assume that $X_1$ is triconnected. By \eqref{equ:delta van} we hence have that $\delta_kX_0=\delta_kX_1=0$ for $k\geq 3$.
In particular the part \eqref{equ:g4} of loop order $g+4$ of the equations above becomes 
\[
\delta_2 X_2 + \delta_0 X_4=0.
\]
Now since $X_0$ is triconnected, $\delta_2 X_0$ is biconnected. Since $\GC_n^{bi}\subset \GC_n$ is a quasi-isomorphism, we may employ Lemma \ref{lem:aux} to find $Z\in \GC_n$ of loop order $g+2$ such that $X_2+\delta_0 Z\in \GC_n^{bi}$.
Then we may replace $X$ by $X+\delta Z$, so that we can assume from the start that $X_2\in \GC_n^{bi}$. This implies by \eqref{equ:delta van} that $\delta_3X_2=0$.

We then take for our desired element $X'$:
\[
X' = X_0 + X_2 + X_4.
\]
This satisfies $\tilde \delta X'\in \GC_n^{\geq g+6\lp}$, with the new equation in loop order 5:
\[
\delta_5 X_0 + \delta_3 X_2 = 0+ 0=0.
\]
The proof that $d_3=0$ can be obtained similarly, by just truncating the above argument.
\end{proof}

\end{document}